\documentclass[12pt]{amsart}
\newtheorem{theorem}{Theorem}[section]
\newtheorem{proposition}[theorem]{Proposition}
\newtheorem{lemma}[theorem]{Lemma}

\theoremstyle{remark}
\newtheorem{remark}[theorem]{Remark}
\theoremstyle{definition}
\newtheorem{definition}[theorem]{Definition}

\newcommand{\C}{\mathbb C}

\newcommand{\Z}{\mathbb Z}

\newcommand{\SL}{\mathit{SL}}
\newcommand{\tr}{\mathit{tr}}
\begin{document}

\title[Reidemeister torsion for the figure-eight knot]
{Reidemeister torsion of a 3-manifold obtained by a Dehn-surgery along the figure-eight knot}

\author{Teruaki Kitano}


\address{Department of Information Systems Science, 
Faculty of Science and Engineering, 
Soka University, 
Tangi-cho 1-236, 
Hachioji, Tokyo 192-8577, Japan}

\email{kitano@soka.ac.jp}

\maketitle

\begin{abstract}
Let $M$ be a 3-manifold obtained by a Dehn-surgery along the figure-eight knot. 
We give a formula of the Reisdemeisiter torsion of $M$ for any $\SL(2;\C)$-irreducible representation. 
It has a rational expression of the trace of the image of the meridian. 
\end{abstract}

\section{Introduction}

Reidemeister torsion is a piecewise linear invariant for manifolds and originally defined by Reidemeister, Franz and de Rham in 1930's. 
In 1980's Johnson developed a theory of the Reidemeister torsion 
from the view point of certain relation to the Casson invariant of a homology 3-sphere. 
He also derived an explicit formula for the Reidemeister torsion of a homology 3-sphere 
obtained by a $\frac1n$-Dehn surgery along any torus knot for $\SL(2;\C)$-irreducible representations. 
We generalized the Johnson's formula for any Seifert fibered space \cite{Kitano94-1} along his studies. 

In this paper, 
we give a formula for 3-manifolds obtained by Dehn surgeries along the figure-eight knot. 
Let $K\subset S^3$ be the figure-eight knot. 
The knot group $\pi_1(S^3\setminus K)$ has the following presentation. 
\[
\pi_1(S^3\setminus K)=\langle x,y\ |\ wx=yw \rangle
\]
where $w=xy^{-1}x^{-1}y$, $l=w^{-1}\tilde{w}$ and $\tilde{w}=x^{-1}yxy^{-1}$. 
Now $x$ is a meridian and $l$ is a longitude. 

Let $M$ be a 3-manifold obtained by a $\displaystyle\frac pq$-surgery along $K$. 
The fundamental group $\pi_1(M)$ admits a presentation as follows;
\[
\pi_1(M)=\langle x,y\ |\ wx=yw, x^p l^q= 1\rangle.
\]
Let $\rho:\pi_1(M)\rightarrow\SL(2;\C)$ be an irreducible representation. 
Assume the chain complex $C_\ast(M;\C^2_\rho)$ is acyclic.
Then Reidemeister torsion $\tau_\rho(M)=\tau(C_\ast(M;\C^2_\rho))$ is given by the following. 

\begin{theorem}
\[
\tau_\rho(M)
=\frac{2(u-1)}{u^2(u^2-5)}
\]
where $u=\tr(\rho(x))$. 
\end{theorem}

\begin{remark}
We remark the trace $u$ cannot move freely on the complex plane in the above formula. 
The value $u$ depends on the surgery coefficient $p,q$. 
\end{remark}

\flushleft{\bf Acknowledgements}

This research was supported by JSPS KAKENHI 25400101. 
The author thank Michel Boileau, Michael Heusener and Takayuki Morifuji for helpful comments and discussions. 

\section{Definition of Reidemeister torsion}

First let us describe the definition of the Reidemeister torsion for $\SL(2;\C)$-representations. 
Since we do not give details of definitions and known results, 
please see Johnson \cite{Johnson}, Milnor \cite{Milnor61, Milnor62, Milnor66} and Kitano \cite{Kitano94-1,Kitano94-2} for details.

Let $W$ be an $n$-dimensional vector space over $\C$ 
and let ${\bf b}=(b_1,\cdots,b_n)$ and
${\bf c}= (c_1,\cdots,c_n)$ be two bases for $W$. 
Setting $b_i=\sum p_{ji}c_i$, 
we obtain a nonsingular matrix $P=(p_{ij})$ with entries in $\C$. 
Let $[{\bf b}/{\bf c}]$ denote the determinant of $P$.

Suppose
\[
C_*: 0 \rightarrow C_m
\overset{\partial_m}{\rightarrow} 
C_{m-1}
\overset{\partial_{m-1}}{\rightarrow} 
\cdots 
\overset{\partial_{2}}{\rightarrow}
C_1\overset{\partial_1}{\rightarrow} C_0\rightarrow 0
\]
is an acyclic chain complex of finite dimensional vector spaces over $\C$. 
We assume that a preferred basis ${\bf c}_i$ for $C_i$ is given for each $i$. 
Choose some basis ${\bf b}_i$ for $B_i=\mathrm{Im}(\partial_{i+1})$ 
and take a lift of it in $C_{q+1}$, which we denote by $\tilde{\bf b}_i$.
Since $B_i=Z_i=\mathrm{Ker}{\partial_i}$, 
the basis ${\bf b}_i$ can serve as a basis for $Z_i$. 
Furthermore since the sequence
\[
0\rightarrow Z_i
\rightarrow C_i
\rightarrow B_{i-1}\rightarrow 0
\]
is exact, the vectors $({\bf b}_i,\tilde{\bf b}_{i-1}) $ form a basis for $C_i$. 
Here $\tilde{\bf b}_{i-1}$ is a lift of ${\bf b}_{i-1}$ in $C_i$. 
It is easily shown that
$[{\bf b}_i,\tilde{\bf b}_{i-1}/{\bf c}_i]$ does not depend on the choice 
of a lift $\tilde{\bf b}_{i-1}$. 
Hence we can simply denote
it by $[{\bf b}_i, {\bf b}_{i-1}/{\bf c}_i]$.

\begin{definition}
The torsion $\tau(C_*)$ is given by the alternating product 
\[
\prod_{i=0}^m[{\bf b}_i, {\bf b}_{i-1} /{\bf c}_i]^{(-1)^{i+1}}.
\]
\end{definition}

\begin{remark}
It is easy to see that $\tau(C_\ast)$ does not depend on the choices of the bases $\{{\bf b}_0,\cdots,{\bf b}_m\}$. 
\end{remark}

Now we apply this torsion invariant of chain complexes to the following geometric situations. 
Let $M$ be a finite CW-complex and $\tilde M$ a universal covering of $M$. 
The fundamental group $\pi_1(M)$ acts on $\tilde M$ as deck transformations. 
Then the chain complex $C_*(\tilde{M};\Z)$ has the structure of a chain complex of free $\Z[\pi_1(M)]$-modules. 
We denote the 2-dimensional vector space $\C^2$ by $V$. 
Using a representation $\rho:\pi_1(M)\rightarrow \SL(2; \C)$, 
$V$ has the structure of a $\Z[\pi_1(M)]$-module.  
Then we denote it by $V_\rho$ and 
define the chain complex $C_*(M; V_\rho)$ 
by $C_*(\tilde{M}; \Z)\otimes_{\Z[\pi_1(M)]} V_\rho$. 
Here we choose a preferred basis
\[
\{\tilde{u}_1\otimes {\bf e}_1, \tilde{u}_1\otimes {\bf e}_2, \cdots,\tilde{u}_k\otimes{\bf e}_1, \tilde{u}_k\otimes{\bf e}_2
\}
\]
of $C_q(M; V_\rho)$ where $\{{\bf e}_1 , {\bf e_2}\}$ is a canonical basis of $V=\C^2$ 
and 
$u_1,\cdots,u_k$ are the $q$-cells giving the preferred basis of $C_q(M; \Z)$.

We suppose that all homology groups $H_*(M; V_\rho)$ are vanishing. 
In this case we call $\rho$ an acyclic representation. 

\begin{definition}
Let $\rho:\pi_1(M)\rightarrow \SL(2; \C)$ be an acyclic representation. 
Then the Reidemeister torsion $\tau_\rho(M)$ is defined 
to be the torsion $\tau(C_*(M; V_\rho))$. 
\end{definition}

\begin{remark}
\noindent
\begin{enumerate}
\item
We define the $\tau_\rho(M)=0$ for a non-acyclic representation $\rho$.

\item
The Reidemeister torsion $\tau_\rho(M)$ depends on several choices. 
However it is well known that the Reidemeister torsion is a piecewise linear invariant. 
See Johnson \cite{Johnson} and Milnor \cite{Milnor61, Milnor62, Milnor66}.
\end{enumerate}
\end{remark}

Here we recall the Reidemeister torsion of the torus and the solid torus.  

\begin{proposition}
Let $\rho:\pi_1(T^2)\rightarrow\SL(2; \C)$ be a representation. 
\begin{enumerate}
\item
This representation $\rho$ is an acyclic representation if and only if there exists an element $z\in\pi_1(T^2)$ such that $\mathrm{tr}(\rho(z))\neq 2$. 
\item
If $\rho$ is acyclic, then it holds $\tau_\rho(T^2)= 1$.
\end{enumerate}
\end{proposition}

Next we consider the solid torus $S^1\times D^2$ with $\pi_1 (S^1\times D^2)\cong\Z$ generated by $x$.

\begin{proposition}
Let $\pi_1(S^1\times D^2)\rightarrow \SL(2;\C)$ be a representation. 
Then it holds  
\[
\begin{split}
\tau(S^1\times D^2; V_\rho)
&=\frac{1}{\mathrm{det}(\rho(l)-E)}\\
&=\frac{1}{2-\mathrm{tr}(\rho(l))}
\end{split}
\]
for a generator $l\in\pi_1(S^1\times D^2)\cong\Z$. 
Here $E$ is the identity matrix in $\SL(2;\C)$. 
\end{proposition}

From here we assume $M$ is a compact 3-manifold 
with an acyclic representation $\rho:\pi_1(M)\rightarrow \SL(2; \C)$.  
Here we take a torus decomposition of $M=A\cup_{T^2} B$. 
For simplicity, we write the same symbol $\rho$ 
for a restricted representation to subgroups 
$\pi_1(A),\pi_1(B)$ and $\pi_1(T^2)$ of $\pi_1(M)$. 

By this torus decomposition, 
we have the following exact sequence:
\[
0\rightarrow
C_\ast(T^2;V_\rho)
\rightarrow
C_\ast(A;V_\rho)\oplus C_\ast(B;V_\rho)
\rightarrow
C_\ast(M;V_\rho)
\rightarrow 0. 
\]

\begin{proposition}
Let $\rho:\pi_1(M)\rightarrow \SL(2;\C)$ be a representation 
which restricted to $\pi_1(T^2)$ is acyclic. 
Then 
$H_\ast(M;V_\rho)=0$ if and only if $H_\ast(A;V_{\rho})=H_\ast(B;V_{\rho})=0$. 
In this case it holds 
\[
\tau_\rho(M)=\tau_{\rho}(A)\tau_{\rho}(B).
\]
\end{proposition}

We apply this proposition to any 3-manifold obtained by Dehn-surgery along a knot. 
Now let $M$ be a closed 3-manifold 
obtained by a $\frac{p}{q}$-surgery along the figure eight knot $K$. 
Under the presentation  
\[
\pi_1(E(K))=\langle x,y\ |\ wx=yw \rangle
\]
where $w=xy^{-1}x^{-1}y$, $l=w^{-1}\tilde{w}$ and $\tilde{w}=x^{-1}yxy^{-1}$, 
$x$ is a meridian and $l=w^{-1}\tilde{w}$ is a longitude. 

We take an open tubular neighborhood $N(K)$ of $K$ and its knot exterior $E(K)=S^3\setminus {N}(K)$. 
We denote its closure of ${N(K)}$ by $\bar{N}$ which is homeomorphic to $S^{1}\times D^{2}$.
Since this 3-manifold $M$ is obtained by Dehn-surgery along $K$, 
we have a torus decomposition 
\[
M=E(K)\cup \bar{N}. 
\]

Let $\rho:\pi_1(E(K))=\pi_1(S^3\setminus K)\rightarrow\SL(2;\C)$ be a representation 
which extends to $\pi_1(M)$. 
In this case it holds the following.

\begin{proposition}
If $\rho$ is acyclic on $\pi_1(T^2)$, then 
$\tau_\rho(M)
=\tau_\rho(E(K))\tau_\rho(\bar{N})$. 
Further if all chain comeplexes are acyclic, then 
\[\tau_\rho(M)
=\frac{\tau_\rho(E(K))}{2-\mathrm{tr}(\rho(l))}.
\]
\end{proposition}

\section{Main result}

Recall the following lemma, which is the fundamental way 
to study $\SL(2;\C)$-representations of a 2-bridge knot. 
Please see \cite{Riley84} as a reference. 

\begin{lemma}
Let $X,Y\in\SL(2,\C)$. If $X$ and $Y$ are conjugate and $XY\neq YX$, 
then there exists $P\in\SL(2;\C)$ s.t. 
\[
PXP^{-1}=\begin{pmatrix} s & 1\\0 & 1/s\end{pmatrix},\ 
PYP^{-1}=\begin{pmatrix} s & 0\\-t & 1/s\end{pmatrix}.
\]
\end{lemma}

We apply this lemma to irreducible representations of $\pi_1(E(K))$. 
For any irreducible representation $\rho$, 
we may assume that its representative of this conjugacy class is given by
\[
\rho_{s,t}:\pi_1(E(K))\rightarrow \SL(2;\C)\ (s,t\in \C\setminus\{0\})\]
where
\[
\rho_{s,t}(x)=\begin{pmatrix} s & 1\\0 & 1/s\end{pmatrix},
\rho_{s,t}(y)=\begin{pmatrix} s & 0\\-t& 1/s\end{pmatrix}
\]
Simply we write $\rho$ to $\rho_{s,t}$ for some $s,t$. 
We compute the matrix 
\[
R=\rho(w)\rho(x)-\rho(y)\rho(w)=(R_{ij})
\]
to get the defining equations of the space of the conjugacy classes of the irreducible representations.

\begin{itemize}
\item
$R_{11}=0$,
\item
$R_{12}=3-\frac{1}{s^2}-s^2+3 t-\frac{t}{s^2}-s^2 t+t^2$,
\item
$R_{21}=3 t-\frac{t}{s^2}-s^2 t+3 t^2-\frac{t^2}{s^2}-s^2 t^2-t^3=tR_{12}$,
\item
$R_{22}=0$.
\end{itemize}

Hence $R_{12}=0$ is the equation of the space of the conjugacy classes of the irreducible representations. 

This equation 
\[
3-\frac{1}{s^2}-s^2+3 t-\frac{t}{s^2}-s^2 t+t^2=0
\]
can be solved in $t$ as 
\[
t=\frac{1-3 s^2+s^4\pm \sqrt{1-2 s^2-s^4-2 s^6+s^8}}{2 s^2}.
\]

Here it can be seen that $L=\rho(l)=(l_{ij})$ is given by the followings:
\begin{lemma}
\[
\begin{split}
l_{11}
&=
1-\frac{t}{s^2}+s^2 t-t^2+\frac{t^2}{s^4}-\frac{t^2}{s^2}+s^2 t^2-t^3-\frac{t^3}{s^2}\\
l_{12}
&=
\frac{t}{s^3}+s^3 t-\frac{t^2}{s}-st^2
\\
l_{21}
&=
\frac{t^2}{s^3}-\frac{2t^2}{s}-2 st^2+s^3 t^{2}+\frac{t^3}{s^3}-\frac{2t^3}{s}-2 st^3+s^3 t^3-\frac{t^4}{s}-st^4\\
l_{22}
&=
1+\frac{t}{s^2}-s^{2}t-t^2+\frac{t^2}{s^2}-s^2 t^2+s^4t^2-t^3-s^2 t^3
\end{split}
\]
\end{lemma}

Here we get the trace of direct computation. 
\[
\mathrm{tr}(\rho(l))
=2-2t^2+\frac{t^2}{s^4}+s^4 t^2-2 t^3-\frac{t^3}{s^2}-s^2 t^3
\]

It is easy to see that $\mathrm{tr}(\rho(l))\neq 2$ if $u=s+\frac{1}{s}=2$. 
Hence there exists an element $z\in\pi_1(T^2)$ s.t. $\mathrm{tr}(\rho(z)\neq 2$. 
This means $\rho$ is always acyclic on $T^2$. 
Now we have 
\[
\tau_\rho(M)
=\tau_\rho(E(K))\tau_\rho(\bar{N}).\]

Here we obtain the Reidemeister torsion of $E(K)$ as follows. 
See \cite{Kitano94-2} for precise computation. 

\begin{proposition}
\[
\tau_\rho(E(K))=-2(u-1)
\]
where $u=s+\frac{1}{s}$.
\end{proposition}

By substituting 
\[
t=\frac{1-3 s^2+s^4\pm \sqrt{1-2 s^2-s^4-2 s^6+s^8}}{2 s^2}
\]
in $\mathrm{tr}(\rho(l)$, 
we get the following proposition. 

\begin{proposition}
\[
\tau_\rho(\bar{N})
=-\frac{1}{u^2(u^2-5)}.
\]
\end{proposition}

Therefore we obtain the following formula:
\[
\begin{split}
\tau_\rho(M)
&=\tau_\rho(E(K))\tau_\rho(\bar{N})\\
&=(-2(u-1))\left(-\frac{1}{u^2(u^2-5)}\right)\\
&=\frac{2(u-1)}{u^2(u^2-5)}.
\end{split}
\]

\begin{remark}
The representations for $u^2-5=0$ are degenerate into reducible representation from irreducible representations. 
\end{remark}


\end{document}